\documentclass[11pt]{amsart}      

\usepackage{amssymb}      
\usepackage{amsmath}      
\usepackage{amsthm}      
\usepackage{amsbsy}       
\usepackage{amscd} 
\usepackage{latexsym} 
\usepackage[dvips]{graphicx}

\usepackage{bm}    

\usepackage[
	hypertexnames=false,
	hyperindex,
	pagebackref,  
	bookmarks=false,
	colorlinks,
	linkcolor=blue,
	citecolor=red,
	urlcolor=red,
]{hyperref}

\theoremstyle{plain}      
\newtheorem{theorem}{Theorem}[section]      

\theoremstyle{remark}      
\newtheorem{remark}{Remark}[section]

\newcommand{\Q}{{\mathbb{Q}}}        
\newcommand{\Z}{{\mathbb{Z}}}   
 \newcommand{\N}{{\mathbb{N}}}  
\newcommand{\C}{{\mathbb{C}}}      
\newcommand{\R}{{\mathbb{R}}}

 \newcommand{\M}{{\mathcal{M}}}

 \newcommand{\Hom}{{\rm{Hom}}}  
\newcommand{\G}{{\Gamma}}

\newcommand{\ldual}[1]{#1^{*}}

\newcommand{\A}{\mathfrak A}

  \newcommand{\ro}{{\widetilde{\rho}}}

\begin{document}


\title[Mapping class groups]{On mapping class groups and their TQFT representations    
}     

 \author{Louis Funar} 
 \address{ Univ. Grenoble Alpes, CNRS, Institut Fourier, 38000 Grenoble, France}
 \email{louis.funar@univ-grenoble-alpes.fr}

\maketitle

\begin{abstract}
The aim of this paper is to survey some aspects of mapping class groups with focus on their finite dimensional representations arising in topological quantum field theory.

\vspace{0.1cm}
 
\end{abstract}

\section{Introduction}
Mapping class groups form a well-established subject in mathematics at the interface of low dimensional topology, algebraic topology,  
dynamics and algebraic geometry. It is a rich and active domain of mathematics, as it can be inferred from the collection of articles in \cite{Farb}, where a large number of open questions are discussed. We refer to  \cite{FaMa} for a thoroughful introduction to the subject. Recall that the mapping class group of a  manifold is the group of isotopy 
classes of its homeomorphisms. The same terminology is also used for the group of isotopy classes of diffeomorphisms of 
smooth manifolds. While these definitions agree in the case of compact surfaces, which is the main concern here,  
this is not necessarily the case for higher dimensional manifolds.  Mapping class groups of low dimensional compact manifolds 
are countable groups which have been extensively studied in relation with Teichm\"uller spaces and hyperbolic geometry.  Our aim is to report on some aspects of this subject which were strongly influenced by the seminal work of Witten (see \cite{Witten}) in quantum field theory and are today part of the research field called quantum topology. 

The objectives of this chapter are as follows: 
\begin{enumerate}
\item introduce TQFTs and universal pairings
\item analyse the positivity of universal pairings and manifold detection by TQFTs
\item introduce modular tensor categories 
\item state some properties of modular representations in genus one, e.g. the congruence property
\item discuss the modular representations in higher genera for the skein TQFT
\end{enumerate}

Topological quantum field theories, abbreviated TQFTs in the sequel, have their origins in the paper \cite{Atiyah} of Atiyah who formulated a number of axioms inspired from cobordism theory,  Jones' discovery of a new link invariant (\cite{Jones}), the work of Witten in quantum field theory and 
Segal's mathematical description of the conformal field theory in dimension 2 (\cite{Segal,Segal2}).  Specifically, for  a fixed $n$ 
we consider the category whose objects are homeomorphisms classes of $(n-1)$-dimensional closed manifolds and 
morphisms correspond to homeomorphism classes of cobordisms between them. Other versions of this category arise when  
manifolds are endowed with some extra structure and the corresponding homeomorphisms need to respect it. For instance,  
we can restrict to smooth manifolds and their diffeomorphisms, to oriented manifolds and orientation-preserving 
homeomorphisms and so on.  
Then an $n$-dimensional  (topological) TQFT 
is a functor from  this  category  into the category of finite dimensional vector spaces and linear maps between them. 
Namely we have a vector space $V(S)$ associated to each closed $(n-1)$-manifold $S$ and a linear map 
$Z(M)\in \Hom(V(\partial_+ M),V(\partial_-M))$ associated to each cobordism $M$ between 
the manifolds  $\partial_+ M$ and  $\partial_-M$ (in this order)  which only depends on the homeomorphism type 
of the cobordism and the boundary identifications. The main difference with the classical functors encountered in algebraic topology is the monoidal axiom below:  
\[ V(M_1\sqcup M_2) = V(M_1) \otimes V(M_2)\]
where $\sqcup$ denotes the disjoint union. As alluded above we can consider $n$-dimensional smooth TQFT, by restricting 
to smooth manifolds and cobordisms, while asking $Z(M)$ to be an invariant of the diffeomorphism type of the cobordism $M$. 
One usually adds the axiom $V(\emptyset)=\C$, such that the TQFT functor  associates a linear form $Z(M)\in V(\partial M)^*$ 
when $\partial_+M=\partial M$ and a homomorphism  $Z(M)\in \Hom(\C,V(\partial M))$, when  $\partial_+M=\emptyset$. In the later case 
by notation abuse we denote  the image of $Z(M)(1)$ of $1\in \C$ as $Z(M)\in V(\partial M)$. Note that for a closed $n$-manifold, viewed as a cobordism between two empty sets the functor associates a complex number $Z(M)\in \C$, which is a topological invariant of the manifold. If the theory is only defined for a certain class of closed $(n-1)$-manifolds, we say that we have a {\em restricted} TQFT.

Soon after Atiyah's advertisement of topological quantum field theories Witten provided the first example of a meaningful TQFT in \cite{Witten} in dimension 3, known as the  Chern-Simons theory with gauge group $SU(2)$.  The invariant that this TQFT associates to the link complement in the 3-sphere is essentially the coloured Jones polynomial of the link. The mathematically rigorous construction of the $SU(2)$-Chern-Simons theory and its generalisations throughout quantum groups was obtained in the work \cite{RT} of Reshetikhin and Turaev. 

TQFTs in dimension 2 correspond to commutative Frobenius algebras. Constructions of TQFT in dimension 3 rely on some algebraic data called modular tensor categories  (see \cite{Turaev}), which are defined in section \ref{section3},  the main examples coming from quantum groups.   In contrast with the dimension 3 where TQFTs are abundant, TQFTs in higher dimensions are scarce.  
Although it is still unknown what about the topology they detect, TQFTs  in dimension 3 are an useful tool to understand the 
the representation theory of surface automorphisms. 

The simplest  $n$-dimensional cobordisms are the mapping cylinders $C(\phi)$ of homeomorphisms $\phi:M\to M$ 
of $(n-1)$-dimensional manifolds $M$. Strictly speaking $C(\phi)$ is topologically trivial, as it is homeomorphic to 
$M\times [0,1]$, however the boundary identifications with $M$ encode the twist by $\phi$. Therefore any TQFT $Z$ will provide 
a representation of the group of homeomorphisms of $M$ into the group $GL(V(M))$ of automorphisms of the vector space $V(M)$. 
As we required the theory to be topologically invariant this representation factors through 
the mapping class group $\M(M)$ of $M$. In particular a TQFT in dimension $n$ subsumes a family of mapping class group representations  $\rho_M$, for all closed $(n-1)$-manifolds $M$. 
The dimension $n=3$ TQFT will therefore 
provide a large supply of surface mapping class group representations.  The goal of this note is to analyse some of the 
known results about these representations and to promote them for further study. 

We refer to \cite{Atiyah} for a complete set of axioms underlying a TQFT. Let us restrict to the category of 
oriented manifolds and orientation-preserving homeomorphisms  and $\overline{N}$ denote the manifold $N$ with the opposite orientation. Then a TQFT is Hermitian if vector spaces are endowed with non-degenerate Hermitian forms and 
$Z(\overline{M})=\overline{Z(M)}$, where we identify $V(\overline{M})$ with the dual Hermitian space $V(M)^*$. Among Hermitian TQFT 
those having a physical meaning are the unitary ones, namely those for which the  Hermitian forms are positive definite. 

A key assumption on the TQFTs is that all vector spaces $V(M)$ are {\em finite dimensional}. Note that, depending on the type of TQFT one considers, this is not always an assumption, but can arise automatically. 
Then the value of the invariant $Z(N)$ of a closed $n$-manifold which fibers over the circle $S^1$ is given by: 
\[ Z(N)= {\rm Tr}(\rho_M(\phi)),\]
where $M$ is the fiber, $\phi\in \M(M)$ is the monodromy of the circle fibration and $\rho_M$ is the mapping class group representation above. In particular
\[Z(M\times S^1)= \dim V(M)\]
By abuse of language we will refer to  {\em infinite dimensional} TQFT, if we drop the finiteness assumption below.

\section{Universal pairings and TQFTs}
\subsection{Universal pairings} The TQFT setting can be realised at universal level by some abstract nonsense constructions. 
Assume for now that we work with the category of smooth manifolds. 
For each  closed oriented $(n-1)$-manifold $M$ let $U_M$ be the complex vector space spanned by the set of classes of 
compact connected $n$-manifolds $N$ with a marking diffeomorphism $\partial N\to M$, up to diffeomorphisms commuting with the 
boundary markings. 
If $M=\emptyset$ we set $\A=U_{\emptyset}$. Note that $\A$ has a natural algebra structure induced by the 
connect-sum operation for closed manifolds. 

The {\em universal pairing associated to $M$} is the complex sesquilinear map
\[\langle , \rangle_M: U_M\times U_M \to \A\]
which extends by sesquilinearity the map sending two marked $n$-manifolds $N_1,N_2$ with boundary $M$ to 
the result of their gluing along the boundary $N_1\cup_M \overline{N_2}$.  

\subsection{Functors from topological invariants} This pairing provides already a recipe for constructing some functors similar to TQFTs out of topological invariants. 
A complex valued topological invariant $I$ of closed orientable $n$-manifolds which is sensitive to the orientation, 
i.e. such that $I(\overline{N})= \overline{I(N)}$,  extends by linearity to  $I: \A\to \C$.  
For connected $M$, the quotient vector space:
\[ V_I(M)= U_M/\ker(I\circ \langle \;,\;\rangle_{M})\]
by the left kernel of the Hermitian form  $I\circ \langle \;,\;\rangle_{M}$ inherits a non-degenerate Hermitian form. 
We extend the definition of $V_I(M)$ to non connected manifolds by means of the monoidal axiom. Then for every $n$-dimensional cobordism 
$N$ between $M_1$ and $M_2$  we associate the linear map $Z_I(N)\in \Hom(V_I(M_1),V_I(M_2))$ with the property that for every 
$N_1$ a $n$-manifold with boundary $M_1$ we have:
\[ Z_I(N)([N_1])= [N\cup_{M_1} N_1]\in V_I(M_2)\]
Then $Z_I$ is a functor from the category of $(m-1)$-dimensional closed oriented smooth manifolds and cobordisms between them into the 
category of vector spaces endowed with non-degenerate Hermitian forms.  One might need extra assumptions on $I$ in order that 
$V_I(M_1\sqcup M_2)=V_I(M_1)\otimes V_I(M_2)$, for all manifolds $M_i$. However, the functoriality of $Z_I$ provides a family of mapping class groups representations $\rho_{I,M}$ of $\M(M)$ into the group of Hermitian automorphisms of $V_I(M)$. 
Although the functor $Z_I$ is close to a TQFT, most often the vector spaces  $V_I(M)$ are infinite dimensional unless $I$ is of a very particular form. 
In the reverse direction a family of finite dimensional surface mapping class group representations satisfying some compatibility conditions can be extended to a TQFT in dimension 3 (see \cite{F95,Juhasz}).

\subsection{Positivity of universal pairings and manifold detection}
A natural question is whether the universal pairings are {\em positive definite}, namely if 
$\langle v, v\rangle_M=0$ only if $v=0\in U_M$. 
This is directly related to the ability of unitary TQFTs to distinguish among 
diffeomorphism types of manifolds. Indeed, if the universal pairing is not positive definite and $v\in V_M$ is a vector such that 
$\langle v, v\rangle_M=0$, then no unitary TQFT in dimension $n$ can distinguish the linear combination  of manifolds $v$ 
from the zero vector. By now we have a complete answer to the positivity question.

It is  proved in \cite{FKNSWZ} that universal pairings are positive definite in dimension $n\leq 2$.  
This is not anymore true in dimension $4$, according to the same authors:

\begin{theorem}[\cite{FKNSWZ}]
There are homology 3-spheres $M$ for which the universal pairing $\langle \;, \;\rangle_M$ is not positive definite.  
Moreover, unitary  restricted TQFT on homology 3-spheres cannot distinguish among simply connected closed smooth $4$-manifolds which are homotopy equivalent.  
Furthermore unitary TQFT cannot distinguish among smoothly s-cobordant 4-manifolds. 
\end{theorem}

The failure of TQFTs to distinguish between smooth 4-manifolds was further extended to all higher dimensional manifolds by 
Kreck and Teichner in \cite{KT}: 

\begin{theorem}[\cite{KT}]
When $M=S^{4k-1}$, $k\geq 2$, or $M=S^3\times S^m$, $m\geq 1$, then the universal pairing $\langle \;, \;\rangle_M$ is not positive definite and moreover there exist non-diffeomorphic orientable manifolds $N_1,N_2$ with boundary $M$, which can be chosen simply connected when  $\dim M\neq 4$,   such that: 
\[\langle N_1-N_2, N_1-N_2\rangle_M=0 \] 
\end{theorem}

On the other hand, Barden (\cite{Barden} has shown that simply connected smooth 5-manifolds are classified by their second homology and  the $i$-invariant which is the maximal $r$ for which the second Stiefel-Whitney class $w_2(M)\in H^2(M;\Z/2\Z)$ admits a lift to a class in  $H^2(M;\Z/2^r\Z)$. Krek and Teichner showed that vanishing of cohomology operations can be detected by means of specific TQFTs, the so-called called Dijkgraaf-Witten invariants coming from spaces with finite total homotopy (see \cite{KT}) and consequently we have that: 
\begin{theorem}[\cite{KT}]
Unitary TQFTs detect the homology and the $i$-invariant and in particular  they distinguish simply connected closed oriented 5-manifolds. Therefore, the universal pairing $\langle \;, \;\rangle_{S^4}$ is positive when restricted to  simply connected 5-manifolds. 
\end{theorem}

The remaining dimension 3 seems the most challenging case. 
Eventually Calegari, Freedman and Walker gave a positive answer in dimension 3, in full generality: 

\begin{theorem}[\cite{CFW}]
For all closed oriented surfaces $S$ the universal pairing $\langle \;, \;\rangle_S$ is positive definite. 
\end{theorem}

By carefully choosing linear forms $I:\A\to \C$ separating the 3-manifold diffeomorphism types we 
derive some (possibly infinite dimensional) TQFTs $Z_I$.  One might reasonably hope then that TQFTs 
can distinguish 3-manifolds. However, this is not so, because:
\begin{theorem}[\cite{F13}]\label{ideal}
For every $n$ there exist $n$ torus bundles over the circle which cannot be pairwise distinguished by any 
TQFT obtained from a modular tensor category.  
\end{theorem}

Actually TQFTs being finite dimensional, they only carry topological information which can be recognised at some profinite level. By analogy with the Reidemeister-Singer description of closed 3-manifolds we may define profinite closed 3-manifolds to be stable equivalence classes of double cosets of some profinite completion of mapping class groups with respect to the closure of handlebody groups. We can for instance consider the profinite completion induced by all finite quotients of mapping class groups arising from TQFT representations. 
Then two closed 3-manifolds are not distinguished by TQFTs if they become equivalent as profinite 3-manifolds in this sense. 
This means that they have Heegaard decompositions into two handlebodies, whose gluing maps have the same image 
into all finite quotients of the mapping class group, up to left and right multiplications by elements of the handlebody group.

The Dijkgraaf-Witten invariants  (see the next section) are particular topological invariants coming from TQFTs associated to finite groups.  In the case of closed 3-manifolds they already  determine the profinite completion of the fundamental group. 
Whether closed hyperbolic 3-manifolds are distinguished by their Dijkgraaf-Witten invariants amounts to ask when Kleinian groups are profinitely rigid. A recent result of Yi Liu (\cite{Liu}) shows that there are only finitely many such groups having the same profinite completion.

From now on we will stick to TQFTs for oriented 3-manifolds and the associated surface mapping class group representations. 
We shall see in the next section that this strongly constrains  mapping class group 
representations.

\section{Genus one TQFT representations and the Congruence Property}\label{section3}

\subsection{Fusion  categories}  
We follow the definitions from \cite{BK}. 
A {\em rigid} monoidal category is 
a  monoidal category ${\mathcal C}$ 
with unit object ${\mathbf 1}$ such that to each
object $X\in \mathrm{Ob}({\mathcal C})$ there are associated a dual
object $X^*\in \mathrm{Ob}({\mathcal C})$ and four morphisms
\[ \mathrm{ev}_X \colon X^*\otimes X \to{\mathbf 1},  \qquad \mathrm{coev}_X\colon {\mathbf 1}  \to X \otimes X^*,  \qquad 
  \widetilde{\mathrm{ev}}_X \colon X\otimes X^* \to{\mathbf 1}, \qquad   
\widetilde{\mathrm{coev}}_X\colon {\mathbf 1}  \to X^* \otimes X  \]
such that, for every $X\in
\mathrm{Ob}({\mathcal C})$, the pair $(\mathrm{ev}_X,\mathrm{coev}_X)$ is a 
left duality for $X$ and the pair $(\widetilde{\mathrm{ev}}_X,\widetilde{\mathrm{coev}}_X)$ is a right
duality for $X$, namely: 
\begin{equation*}
({\mathbf 1}_X \otimes \mathrm{ev}_X)(\mathrm{coev}_X \otimes {\mathbf 1}_X)={\mathbf 1}_X  \quad  {\rm {and}} \quad (\mathrm{ev}_X \otimes {\mathbf 1}_{X^*})({\mathbf 1}_{X^*} \otimes \mathrm{coev}_X)={\mathbf 1}_{X^*}
\end{equation*}

A {\em pivotal} category is a rigid monoidal category equipped with 
an isomorphism $j$ of monoidal functors between 
identity and $(-)^{**}$, called pivotal structure. Note that the right duality is determined by the left duality and $j$, as follows:
\begin{equation*}
\widetilde{\mathrm{ev}}_X =
\mathrm{coev}_X(\mathbf 1_{X^*}\otimes j_X^{-1}), \: 
\widetilde{\mathrm{coev}}_X= 
(j_X\otimes \mathbf 1_{X^*})\mathrm{ev}_X 
\end{equation*}

Now, for an endomorphism $f $ of an object $X$ of   
a pivotal category ${\mathcal C}$, one defines the
{\em left/right traces} ${\rm tr}_l(f), {\rm tr}_r(f) \in
\mathrm{End}_{\mathcal C}({\mathbf 1})$ by
\begin{equation*}
{\rm tr}_l(f)=\mathrm{ev}_X({\mathbf 1}_{\ldual{X}} \otimes f) \widetilde{\mathrm{coev}_X  \quad {\text
{and}}\quad
 {\rm tr}_r(f)=\widetilde{\mathrm{ev}}_X( f \otimes {\mathbf 1}_{\ldual{X}}) \mathrm{coev}}_X  
\end{equation*}

A {\em spherical category} is a pivotal category whose left and
right traces are equal, i.e.,  ${\rm tr}_l(f)={\rm tr}_r(f)$ for every
endomorphism $f$ of an object. The left and right dimensions of   $X\in \mathrm{Ob} ({\mathcal C})$
are defined by $ \dim_l(X)={\rm tr}_l({\mathbf 1}_X) $ and $
\dim_r(X)={\rm tr}_r({\mathbf 1}_X) $.  
 Then they are 
denoted  ${\rm tr}(f)$ and called the trace of $f$. The
left (and right) dimensions of an object $X$ are denoted  $\dim(X)$
and called the dimension of $X$.

A monoidal {\em ${\mathbb C}$-linear category} is a monoidal category 
${\mathcal C}$ such that its Hom-sets are complex vector spaces and
the composition and monoidal product of morphisms are bilinear. 
An  object $V\in \rm{Ob}(\mathcal C)$ is called {\em simple} 
if the map ${\mathbb C} \to \mathrm{End}_{\mathcal C}({\mathbf 1}), 
a \mapsto a \, {\mathbf 1}_{\mathbf 1}$  is an algebra isomorphism. 
A monoidal  ${\mathbb C}$-linear category is called {\em semi-simple} if every object  is a direct sum 
of finitely many simple objects with 
finite dimensional Hom spaces and $\mathbf 1$ is a 
simple object. 
Now a {\em fusion} category over $\mathbb C$ is a rigid 
semi-simple $\mathbb C$-linear category ${\mathcal C}$ with 
finitely many simple objects. 

Given a finite group $G$ and a 3-cocycle $\omega\in Z^3(G,\C^*)$ we have an associated  spherical fusion category $Vec^{\omega}_G$ of 
$\omega$-twisted $G$-graded vector spaces. Its objects are finite dimensional  $G$-graded vector spaces $V=\bigoplus_{g\in G} V_g$ and homomorphisms are linear maps respecting the grading. The tensor product $\otimes$ is 
$(V\otimes W)_g=\bigoplus_{a,b\in G; ab=g} V_a\otimes W_b$. Let $\underline{g}$ be 1-dimensional vector space concentrated in degree $g$. Then 
the set of isomorphism types of simple objects corresponds to $\{\underline{g}, g\in G\}$, with  unit object $\mathbf 1=\underline{1}$ and dual $\underline{g}^*=\underline{g^{-1}}$.  If we identify $\underline{g}\otimes \underline{h}$ with $\underline{gh}$, then 
${\rm ev}_{\underline{g}}=\omega(g^{-1},g, g^{-1}){\mathbf 1}_{\mathbf 1}$, ${\rm coev}_{\underline{g}}={\mathbf 1}_{\mathbf 1}$ and the associator map 
$(\underline{g}\otimes \underline{h})\otimes \underline{k} \to \underline{g}\otimes (\underline{h}\otimes \underline{k})$ is $\omega(g,h,k) {\mathbf 1}_{\underline{ghk}}$. The pivotal structure is given by 
$j_{\underline{g}}=\omega(g^{-1},g, g^{-1}){\mathbf 1}_{\underline{g}}$. 

A monoidal category ${\mathcal C}$ is {\em braided} if there exist 
natural isomorphisms $c_{V,W}:V\otimes W\to W\otimes V$ for every objects 
$V,W$,  such that for any $U,V,W\in {\rm Ob}(\mathcal C)$ we have: 
\begin{equation*}
c_{U,V\otimes W}=(\mathbf 1_{V}\otimes c_{U,W})(c_{U,V}\otimes 
\mathbf 1_{W}), \: 
c_{U\otimes V,W}=(c_{U,W}\otimes 
\mathbf 1_{V})(\mathbf 1_{U}\otimes c_{V,W})
\end{equation*}

Let now ${\mathcal C}$ be a left rigid braided monoidal category. 
We do not require that $V^{**}=V$. 
A twist of ${\mathcal C}$ is an automorphism $\theta$ of the identity 
functor of ${\mathcal C}$ satisfying 
\begin{equation*}
\theta_{V\otimes W}=c_{W,V}c_{V,W}(\theta_V\otimes\theta_W), \: \rm{ and } 
\:\: \theta_{\mathbf 1}={\mathbf 1}_{\mathbf 1}
\end{equation*}
The twist $\theta$ is a {\em ribbon} structure on $(\mathcal C, c)$ 
if it also satisfies $\theta_V^*=\theta_{V^*}$ for every 
$V\in \rm{Ob}(\mathcal C)$, and 
the (left) duality is compatible with the ribbon and twist structures, namely:
\begin{equation*}
(\theta_V\otimes \mathbf 1_{V^*}){\rm coev}_V=
(\mathbf 1_V\otimes \theta_{V^*}){\rm coev}_V
 \end{equation*}
In this case $(\mathcal C, c, \theta)$ is called a {\em ribbon} category. 
To a ribbon category one associates naturally a pivotal structure 
by using the (canonical) isomorphism $u_X:X\to X^{**}$ given by:
\begin{equation*}
u_X=({\rm ev}_{X^*}\otimes {\mathbf 1}_X)({\mathbf 1}_{X^*}\otimes c^{-1}_{X,X^{**}})({\rm coev}_X\otimes {\mathbf 1}_{X^{**}})
\end{equation*}
and setting $\theta=u^{-1}j$. 
Moreover this pivotal structure $j$ is spherical. 

A  semisimple {\em  modular tensor category} (MTC)  is a ribbon fusion category $({\mathcal C},c,\theta)$, namely a semisimple ribbon category with finitely many isomorphism classes of simple objects $U_i, i\in I$,  such that the matrix $S$
having entries $S_{ij}=tr(c_{U_j,U_i^*}c_{U_i^*,U_j})$ is non-singular, where 
$i,j\in I$. 

This  matrix is called the {\em  $S$-matrix} of the category ${\mathcal C}$. The topological interpretation of the 
entries $S_{ij}$ is the topological invariant associated to the Hopf link coloured by $i$ and $j$. 
Notice that $I$ has induced 
a duality $*$ such that $U_{i^*}=U_i^*$, for any $i\in I$ and 
there exists a label (also called colour) $0\in I$ such that 
$U_0=\mathbf 1$.  
Since the object $U_i$ is simple the twist $\theta_{U_i}$ acts on
$U_i$ as a scalar $\omega_i\in \mathbb C$.  

The centre $Z(\mathcal C)$ of a spherical fusion category is a MTC. For instance the centre of $Vec^{\omega}_G$ 
is the MTC of finite dimensional $D^{\omega}(G)$-modules, where $D^{\omega}(G)$ is the  $\omega$-twisted Drinfeld double 
of the group $G$, namely the quasi-Hopf algebra $\C[G]^*\otimes \C[G]$ with product  twisted by $\omega$: 
\[ (x^*\otimes g)(y^*\otimes h)=\frac{\omega(x,g,h)\omega(g,h,h^{-1}g^{-1}xgh)}{\omega(g,g^{-1}xg,h)} \delta_{x,gyg^{-1}} (x^*\otimes gh)\]  

Now, MTC provide topological invariants of extended 3-manifolds, whose boundary is equipped 
with a subspace of the boundary homology which is Lagrangian with respect to the natural symplectic structure:

\begin{theorem}[\cite{Turaev}]
MTCs provide TQFTs for (extended) 3-manifolds. 
\end{theorem}

Lyubashenko introduced in \cite{Ly}  a more general notion of modular tensor category, not necessarily semisimple. This extends the invariants defined by Hennings (\cite{He}),  Kauffman and Radford (\cite{KR}) and Kuperberg (\cite{Kup}) and works for  the category of modules of any finite dimensional  factorizable ribbon Hopf algebra. Lyubashenko proved that  non semisimple MTCs still provide a family of surface mapping class group representations (\cite{Ly2}, see also \cite{FSS}) and a restricted TQFT, only defined for connected surfaces, which satisfies a weaker version of the monoidal axiom with respect to connected sums (\cite{KL}). A TQFT defined for the category of all surfaces and admissible cobordisms was only recently obtained in \cite{DGGPR} that gives rise 
to the  same mapping class group representations (\cite{DGGPR2}).

\subsection{$SL(2,\Z)$ representations  from MTCs}
The TQFTs associated to MTCs determine therefore representations of mapping class groups. Because of additional structures on surfaces, like the choice of a Lagrangian in the homology, we only derive projective representations, which lift to 
linear representations of central extensions of mapping class groups by finite cyclic groups. In the genus one case these projective representations have (non-unique but finitely many) lifts to linear representations of $SL(2,\Z)$. 

The  {\em projective} representation  in genus one associated to 
the TQFT defined by a modular tensor category ${\mathcal C}$ over $\mathbb C$ has the form 
$\overline{\rho}_{\mathcal C}:SL(2,\Z)\to PGL(K(\mathcal C))$, where $K(\mathcal C)$ is the Grothendieck ring of $\mathcal C$ with $\C$-coefficients, hence additively 
the vector space generated by the isomorphism classes of simple objects of $\mathcal C$, indexed by a finite set $I$.

The group $SL(2,\Z)$ is generated by the matrices 
$s=\left(\begin{array}{cc}
0 & -1 \\
1 & 0 \\
\end{array}
\right)$ and 
$t=\left(\begin{array}{cc}
1 & 1 \\
0 & 1 \\
\end{array}
\right)$. The usual presentation of $SL(2,\Z)$ in terms of the generators 
 $s, t$  is given by the relations 
$(s t)^3=s^2$ and $ s^4=1$.

The matrices entering in the definition of $\overline{\rho}_{\mathcal C}$ 
are the $S$-matrix  defined above and  the $T$-matrix associated to the twist. 
Specifically, $T$ has the entries $T_{ij}=\theta_i\delta_{ij}$, $i,j\in I$.  
Moreover there is also the so-called charge conjugation 
matrix $C$  which is a permutation matrix having entries  
of the form $C_{ij}=\delta_{i\,j^*}$, $i,j\in I$.
 
The Gauss sums of ${\mathcal C}$ are given by 
$p_{\mathcal C}^{\pm}=\sum_{i\in I}\theta_i^{\pm 1}\dim(U_i)^2$ and these are non-zero scalars satisfying: 
\begin{equation*}
p_{\mathcal C}^{+}p_{\mathcal C}^{-}=\sum_{i\in I} \dim(U_i)^2=\dim(\mathcal C)
\end{equation*}

We then have the relations: 
\[
(ST^3)=p_{\mathcal C}^{+} S^2, 
(ST^{-1})^3=p_{\mathcal C}^{-} S^2 C, 
CT=TC, CS=SC, C^2 =1,
S^2=  p_{\mathcal C}^{+}p_{\mathcal C}^{-} C\]

The projective representation 
$\overline{\rho}_{\mathcal C}:SL(2,\Z)\to PGL(K(\mathcal C))$
is defined by 
\begin{equation*}
\overline{\rho}_{\mathcal C}(s)=S, \,\,
\overline{\rho}_{\mathcal C}(t)=T
\end{equation*}

Further one chooses  ${D}\in \mathbb C$ such that 
${D}^2=\dim({\mathcal C})$ and  $\zeta\in \mathbb C$ such that 
$\zeta^3=p_{\mathcal C}^{+}{D}^{-1}$. 
This   enables us to define a lift of $\overline{\rho}_{\mathcal C}$ 
to an ordinary linear representation
\[\rho_{\mathcal C}^{D, \zeta}:SL(2,\Z)\to GL(K(\mathcal C))\] 
by setting 
\begin{equation*}
\rho_{\mathcal C}^{D, \zeta}(s)=D^{-1}S, \,\,
\rho_{\mathcal C}^{D, \zeta}(t)=\zeta^{-1}T
\end{equation*}

These lifts are called the modular representations associated to 
$\mathcal C$.

\subsection{Properties of the $SL(2,\Z)$-representations}
The entries of the  matrices $S$ and $T$ are indexed by $\{0,1,2, \ldots, r-1\}$, where $0$ corresponds to the unit object and $r$ is the number of isomorphism classes of simple objects in the MTC. 
The following collects the most important arithmetic constraints for the entries of the matrices $S$ and $T$ (see e.g. \cite{EGNO}, chap.8 for complete proofs):  
\begin{theorem}
\begin{enumerate}
\item $p_{\mathcal C}^{+}/p_{\mathcal C}^{-}$, $\zeta$ and the  eigenvalues $\theta_i$  of the diagonal matrix $T$ are roots of unity, so $T$ has finite order  (Vafa's theorem).
Moreover, the order  $d$  of $T$ is a divisor of the algebraic integer $(\dim \mathcal C)^{5/2}$ (\cite{E}).
\item $S_{ij}=S_{ji}$,  $T_{00}=S_{00}=1$, while the matrix with entries  $S_{i^*j}$ is the inverse of $S$.
\item $S_{ij}/S_{0j}$, $S_{0i}=\dim(U_i)$, $\dim(\mathcal C)/\dim(U_i)^2$ are cyclotomic integers in  $\Q(\zeta_d)$ (\cite{NS,dBG,CG1}).  
\item $\dim(U_i)$ are totally real cyclotomic integers, i.e. all their conjugates are real (\cite{ENO}).  Thus 
$\dim(U_i)^2$  and $\dim(\mathcal C)$ are totally real positive. 
\item  $S_{i^*j}=\overline{S_{ji}}$  (\cite{ENO}, Prop. 2.12, Rem. 2.13), so that each $SL(2,\Z)$ representation associated to a MTC is  totally unitary, namely all its algebraic conjugates are unitary. 
\item The maps $X_i:\{0,1,\ldots,r-1\}\to \C$ defined by $X_i(j)=S_{ij}/S_{0j}$ span a ring of $C$-valued functions on the set with $r$ elements in which they form a basis.

\end{enumerate}
 
\end{theorem}

Note that in a finite rigid braided tensor category (not necessarily semisimple) twists are quasiunipotent but not necessarily of finite order.

It is presently unknown whether all higher genera mapping class group representations associated to                                                                                                                              MTCs are Hermitian/unitary. However,  it is known that they cannot be totally unitary, as they might be infinite.

Among the finite index subgroups of $SL(2,\Z)$ a prominent place is taken by the congruence subgroups, which have a strong arithmetic flavour. Recall that a {\em congruence} subgroup of $SL(2,\Z)$ is a finite index subgroup containing 
the kernel of the  homomorphism $SL(2,\Z)\to SL(2,\Z/m\Z)$, which reduces mod $m$, for some non-zero integer $m$. 
It is well-known that congruence subgroups are scarce among the finite index subgroups of $SL(2,\Z)$. 
For instance, there are infinitely many quotients of the hyperbolic plane by subgroups of $SL(2,\Z)$ having bounded 
genus, although only finitely many of them correspond to quotients by congruence subgroups of $SL(2,\Z)$. 

A strong restriction  is the following congruence property satisfied by all genus one mapping class group  TQFT representations:

\begin{theorem}[\cite{DLN}]\label{congru}
Kernels of genus one modular representations 
$\rho_{\mathcal C}^{\lambda, \zeta}$ associated to MTCs  are congruence subgroups. 
\end{theorem} 

This was first  proved in several cases by  Coste, Gannon and Bantay for rational 
conformal field theories (see \cite{CG,B}) and further by Ng and Schauenburg for the projective representations associated modular tensor categories which are centres of spherical fusion categories (\cite{NS}), by  Peng Xu for conformal field theories derived from 
vertex operator algebras (see \cite{Xu}) and eventually  by 
Dong, Lin and Ng in \cite{DLN} in the greatest generality above. 

Note that the Theorem \ref{ideal} is a consequence of the congruence subgroup Theorem \ref{congru}.

\section{Higher genus representations}
\subsection{Generalities} 
Denote by $\G_g^n$ the mapping class group of the closed oriented surface $\Sigma_g^n$ of genus $g$ with $n$ punctures and 
$P\G_g^n$ the pure mapping class group, namely the subgroup of those classes of homeomorphisms which fix pointwise the punctures. 

Representations arising from rational conformal field theories have been described by 
Moore and Seiberg in \cite{MS}. Consequently images of Dehn twists in $\G_g^n$ 
have finite orders which are uniformly bounded, independently on $g,n$. 
In the case of unitary MTCs, 
this restriction is explained by the following result due to Bridson, Aramayona and Souto:

\begin{theorem}[\cite{AS,Br}]
If $g\geq 3$ then images of Dehn twists under any finite dimensional representation of $\G_g^n$ into a compact Lie group have finite order, depending only on the Lie group. 
\end{theorem}

\subsection{Dijkgraaf-Witten representations}
One family of MTC which is rather well-understood is the category $Mod-D^{\omega}(G)$ of finite dimensional 
$D^{\omega}(G)$ -modules. The corresponding TQFT is called the twisted Dijkgraaf-Witten TQFT. 

\begin{theorem}[\cite{Gustafson}]
Dijkgraaf-Witten representations of $P\G_g^n$ have finite images. 
\end{theorem}
This improves on previous work \cite{ERW} settling  the case of braid groups and 
\cite{FF} treating the untwisted case. 
A conjecture of Rowell claims that all representations associated to a MTC have finite images if and only if the MTC is weakly integral, namely its  Frobenius-Perron dimension is integral. Recall that for an object $X$ of a fusion category its Frobenius Perron dimension FP$\dim(X)$ is the largest real eigenvalue of the matrix describing the multiplication by $X$ in the Grothendieck ring $K(\mathcal C)$ and FP$\dim(\mathcal C)$ is the sum of FP$\dim(X)^2$ over the set of isomorphisms types of simple objects $X$. One knows that FP$\dim:K(\mathcal C)\to \R$ is an algebra homomorphism (\cite{ENO}) and $\dim=$FP$\dim$ for unitary MTCs. The integrality of  FP$\dim(X)$, for all $X$ is equivalent 
to the fact that we have a category of representations of a  finite dimensional quasi-Hopf algebra.

The mapping class group of a punctured torus is still isomorphic to $SL(2,\Z)$. However the punctured torus Dijkgraaf-Witten representations are not necessarily of congruence type,  as can be seen from the example given in (\cite{FF}, Rem. 3.15). 
It is shown in \cite{FF} that for $g\geq 3$  the image of the Torelli subgroup by the untwisted Dijkgraaf-Witten representations 
is nontrivial iff the group $G$ is nonabelian. Recall that the Torelli subgroup consists of those mapping classes which act trivially on the homology of the surface. 
  
\subsection{K\"ahler groups and compactifications of moduli spaces of curves}
In genus one mapping class group representations correspond to local systems over modular curves and we expect a similar result in higher genera. Some results are known in genus zero, for instance the Burau representations of braid groups at roots of unity 
have geometric descriptions, see e.g. \cite{DMo,Mcmullen}.  
Let $\G_g^n[p]$ be the (normal) subgroup generated 
by $p$-th powers of Dehn twists. We already saw that unitary representations should factor through the quotients of the form 
 $\G_g^n(p)= \G_g^n/\G_g^n[p]$. 
The groups $\G_{g}\slash \G_{g}[p]$ appeared first in \cite{F}. 
In \cite{AF} the authors proved that they are virtually K\"ahler groups when ${\rm g.c.d.}(p,6)=1$, by considering 
the Deligne-Mumford compactifications of the moduli spaces of curves for some level structures for which they are smooth projective varieties. The following result shows that mapping class group representations arising from MTC can be seen as local systems on complex projective varieties which are closed to the Deligne-Mumford moduli space of stable curves (\cite{DM}), see also \cite{DMa}. 

Let $g,n \in \N$ such that $2-2g-n<0$, $(g,n)\not = (1,1), (2,0)$. 
For every $p$ the quotient $P\G_{g}^n/\G_{g}^n[p]$  
is the fundamental group of a compact K\"ahler orbifold ${\overline{\mathcal{M}_{g}^{n}}}[p]$ compactifying 
the moduli stack of $n$-punctured genus $g$ smooth algebraic curves $\mathcal{M}_{g}^{n}$, whose 
moduli space is the  Deligne-Mumford moduli space of stable curves. Moreover, we have: 

\begin{theorem}[\cite{EF}]\label{kahler}
For odd $p\geq 5$ and  even  $p\ge 10$,  $g,n$ as above and $(g,n,p)\neq (2,0,12)$ the orbifold 
${\overline{\mathcal{M}_{g}^{n}}}[p]$  is uniformizable, hence 
$P\G_{g}^n/\G_g^n[p]$ is a K\"ahler group, actually arising as the fundamental group of a smooth complex projective manifold.
\end{theorem}

\subsection{The skein TQFT}\label{tqft} To define the skein TQFT we consider the category whose  
objects are closed oriented surfaces endowed with 
coloured banded points and morphisms are cobordisms 
decorated by uni-trivalent ribbon graphs (i.e. whose vertices have degrees 1 or 3) compatible with the banded points.  A banded point on a surface is a point with a tangent vector at that point.

The TQFT functor $\mathcal V_p$, 
for $p\geq 3$ and a primitive root of unity $A$ of order $2p$ was 
defined in \cite{BHMV} building on  \cite{Lickorish}, see also \cite{Marche}. 
The vector space  associated by the functor $\mathcal V_p$ to a surface 
is called the {\em space of conformal blocks}.  Let 
$H_g$ be a genus $g$ handlebody with $\partial H_g=\Sigma_g$. 
 Assume given a finite set  of banded 
points on $S_g$. Let $\mathfrak G$ be a uni-trivalent ribbon graph embedded in $H_g$ in such a way that 
$H_g$ retracts onto $\mathfrak G$, its univalent vertices are the banded points and it has no other intersections 
with $\Sigma_g$. 

We fix a natural number $p\geq 3$, called the {\em  level} of the TQFT. 
For the sake of simplicity we only consider odd $p$ here, the construction for even $p$ requires only a few modifications (see \cite{BHMV,Marche}). We define the 
{\em set of colours} in odd level $p$ to be $\mathcal C_p=\{0,2,4,\ldots,p-3\}$. 
An edge colouring of $\mathfrak G$ is called {\em $p$-admissible} if the triangle inequality is satisfied at any trivalent vertex of $\mathfrak G$ and the sum of the three colours around a vertex is bounded by $2(p-2)$. 
Then there exists a basis of the space of conformal blocks  $W_{g, (i_1,i_2,\ldots,i_r)}$  associated to the closed surface $\Sigma_g$ with $r$ banded points coloured by $i_1,i_2,\ldots,i_r\in \mathcal C_p$ which is indexed by the set of all $p$-admissible colouring's of $\mathfrak G$ extending the boundary colouring. Note that 
banded points coloured by $0$ do not contribute. 

Observe that an admissible $p$-colouring of $\mathfrak G$ provides an element of the skein module 
$S_{A}(H_g)$ of the handlebody  with banded boundary points coloured $(i_1,i_2,\ldots,i_r)$, evaluated at the  
primitive $2p$-th root of unity $A$ (see \cite{Lickorish3}). This skein element is obtained by cabling the edges of $\mathfrak G$ by the Jones-Wenzl idempotents prescribed by the colouring and  with prescribed  colours at the banded points.  
We suppose that $H_g$ is embedded in a standard way into the $3$-sphere $S^3$, so that the closure of 
its complement is also a genus $g$ handlebody $\overline{H}_g$.  There is then a 
sesquilinear form: 
\[ \langle \;,\; \rangle: S_{A}(H_g)\times S_{A}(\overline{H}_g)\to \C\]
defined by 
\[ \langle x, y \rangle= \langle x \sqcup y \rangle.\]
Here $x\sqcup y$ is the element of $S_{A}(S^3)$ obtained by the disjoint union of  $x$ and $y$ in 
$H_g\cup\overline{H}_g=S^3$, and $\langle \; \rangle: S_{A}(S^3)\to \C$ is the Kauffman bracket invariant (see also \cite{Roberts1}). 

Eventually the space of conformal blocks $W_{g,(i_1,i_2,\ldots,i_r)}$ is the quotient 
$S_{A}(H_g)/\ker \langle\;,\; \rangle$ by the left kernel of the sesquilinear form above (\cite{Lickorish4}). It follows that $W_{g,(i_1,i_2,\ldots,i_r)}$ is endowed with an induced {\em Hermitian form} $H_{A}$. 

The projections of skein elements associated to the $p$-admissible colouring's of a trivalent graph $\mathfrak G$ as above form an orthogonal basis of $W_{g,(i_1,i_2,\ldots,i_r)}$ with respect to $H_{A}$. It is known (\cite{BHMV}) 
that $H_{A}$ only depends on the $p$-th root of unity $\zeta_p=A^2$ and that in this orthogonal basis the diagonal entries belong to the totally real maximal sub field $\Q(\zeta_p+\overline{\zeta_p})$ (after rescaling). 

Let $\mathfrak G'\subset \mathfrak G$ be a uni-trivalent subgraph  whose degree one vertices are coloured, corresponding to a 
sub-surface $\Sigma'$ of $\Sigma_g$ with coloured boundary. The projections in $W_{g,(i_1,i_2,\ldots,i_r)}$ of skein elements associated to 
the $p$-admissible colouring's of $\mathfrak G'$ form an orthogonal basis of the space of conformal blocks 
associated to the surface $\Sigma'$ with coloured boundary components. 

There is a geometric action of the mapping class groups of the handlebodies $H_g$ and $\overline{H}_g$ respectively
on their skein modules and hence on the space of conformal blocks. Note that the mapping class group of the handlebody 
injects naturally into the mapping class group of the boundary surface and it identifies with the subgroup of 
mapping classes of surface homeomorphisms which extend to the handlebody. 
Moreover, these actions extend to a projective  
action $\rho_{g, p, (i_1,\ldots,i_r),A}$ of $\Gamma_g^r$ on $W_{g,(i_1,i_2,\ldots,i_r)}$ respecting the Hermitian form $H_{\zeta_p}=H_A$ (see \cite{BHMV,Roberts1}). When referring to $\rho_{g, p, (i_1,\ldots,i_r),A}$ the subscript 
specifying the genus $g$ will most often be dropped.
Notice that the mapping class group of an essential (i.e. without annuli or disks complements) 
subsurface $\Sigma'\subset \Sigma_g$ is a subgroup of $\Gamma_g$ which preserves the subspace of conformal blocs 
associated to $\Sigma'$ with coloured boundary.  It is worthy to note that  $\rho_{p,(i_1,\ldots,i_r), A}$ only depends on 
$\zeta_p=A^2$, so we can unambiguously shift the notation for this representation to $\rho_{p,(i_1,\ldots,i_r),\zeta_p}$. 

There is a central extension $\widetilde{\Gamma_g}$ of $\Gamma_g$ by $\Z$ and a linear representation  
 $\ro_{p, \zeta_p}$ on $W_{g}$ which  resolves the projective ambiguity of $\rho_{p,\zeta_p}$. 
 The largest such central extension has class $12$ times the Euler class (see \cite{Gervais,MRo}), but we can restrict to an index $12$ subgroup of it, called $\widetilde{\Gamma}_1$  in \cite{MRo}. When $g\geq 4$ it is a perfect group which coincides with the universal central extension of $\G_g$. 
We consider a subsurface $\Sigma_{g,r}\subset \Sigma_{g+r}$ whose complement consists of $r$ copies of $\Sigma_{1,1}$. 
Let $\widetilde{\Gamma_g^r}$ be the pull-back of the central extension $\widetilde{\Gamma_g}$ to 
the subgroup $\Gamma_{g,r}\subset \Gamma_{g+r}$. Then $\widetilde{\Gamma_{g,r}}$ is also a central extension, which we denote $\widetilde{\Gamma_{g}^r}$ of $\Gamma_g^r$ by $\Z^{r+1}$. 
From \cite{Gervais,MRo} we derive that $\widetilde{\Gamma_g^r}$ is perfect, when $g\geq 3$ and of 
order $10$, when $g=2$.

Let $p\geq 5$ be odd and $\zeta_p$  a primitive 
$p$-th root of unity. We denote 
by $\ro_{p, \zeta_p, (i_1,i_2,\ldots,i_r)}$ the linear representation of the central extension 
$\widetilde{\Gamma_g^r}$ which acts on the vector space $W_{g,p,(i_1,i_2,\ldots,i_r)}$ associated by the TQFT to the surface with the corresponding coloured banded points (see \cite{Gervais,MRo}).

The functor $\mathcal V_p$ associates to a handlebody $H_g$ the projection of the skein element 
corresponding to the trivial colouring of the trivalent graph $\mathfrak G$ by $0$. The Reshetikhin-Turaev invariant associated 
to a closed 3-manifold is given by pairing the two vectors associated to handlebodies in a Heegaard decomposition 
of some genus $g$ and taking into account the twisting by the gluing mapping class action on $W_g$. This is well-defined 
up to a $p$-th root of unity and can be lifted to a topological invariant of 3-manifolds with a $p_1$-structure. 
The Turaev-Viro invariant is the squared absolute value of the Reshetikhin-Turaev invariant. 

One should notice that the skein TQFT $\mathcal V_p$ is unitary, in the sense that $H_{\zeta_p}$ is 
a positive definite Hermitian form when $\zeta_p=(-1)^p \exp\left(\frac{2\pi i}{p}\right)$, corresponding to  
$A_p=(-1)^{\frac{p-1}{2}}\exp\left(\frac{(p+1)\pi i}{2p}\right)$. 
For the sake of notational simplicity, from now we will drop the subscript $p$ 
in $\zeta_p$, when the order of the root of unity will be clear from the context.  
Note that for a general primitive $p$-th root of unity, 
the isometries of $H_{\zeta}$ form a pseudo-unitary group.

Now, the image $\rho_{p,\zeta}(T_{\gamma})$  of a right hand Dehn twist $T_{\gamma}$  in 
a convenient basis given by a trivalent graph is easy to describe. 
Assume that the simple curve $\gamma$ is the boundary of a small disk intersecting once transversely 
an edge  $e$ of the graph $\mathfrak G$. Consider $v\in W_g$ be a vector of the basis given by colourings of the graph $\mathfrak G$ 
and assume that edge $e$ is labelled by the colour $c(e)\in \mathcal C_p$. Then the action of 
the (canonical)  lift  $\widetilde{T_{\gamma}}$ of the Dehn twist $T_{\gamma}$ in $\widetilde{\Gamma_g}$ 
is given by (see \cite{BHMV}, 5.8) :
\[ \ro_{p,\zeta}(\widetilde{T_{\gamma}}) v =A^{c(e)(c(e)+2)} v\]

\subsection{Unitary groups of spaces of conformal blocks}

For a prime $p\geq 5$  we denote 
by ${\mathcal O}_p$ the  ring of cyclotomic integers 
${\mathcal O}_p=\Z[\zeta_p]$, if   
$p\equiv 3({\rm mod}\: 4)$ and ${\mathcal O}_p=\Z[\zeta_{4p}]$, if  
$p\equiv 1({\rm mod}\:4)$ respectively, where $\zeta_r$ denotes a primitive
$r$-th root of unity.

Let $S\mathbb U_{g,p,\sigma(\zeta)}$ be the 
special unitary group associated to the 
Hermitian form conjugated by $\sigma$, thus corresponding to some 
Galois conjugate root of unity.
Consider the  semi-simple algebraic group $\mathbb G_{g,p}$ obtained 
by the so-called restriction of scalars construction from the 
totally real cyclotomic field $\Q(\zeta_p+\overline{\zeta_p})$ to $\Q$. 
Specifically, the group  $\mathbb G_{g,p}$ is a product 
$\prod_{\sigma\in S(p)}S\mathbb U_{g,p,\sigma(\zeta)}$. Here $S(p)$ stands for 
a set of representatives of the classes of complex 
embeddings $\sigma$ of $\mathcal O_p$ 
modulo complex conjugation, or equivalently the set of places of 
the totally real  cyclotomic field $\Q(\zeta_p+\overline{\zeta_p})$. 
Denote   
by $\ro_p$ and $\rho_p$ the representations  
$\prod_{\sigma\in S(p)} \ro_{p,\sigma(A_p^2)}$ and 
$\prod_{\sigma\in S(p)} \rho_{p,\sigma(A_p^2)}$, respectively. 
Notice that the real Lie group $\mathbb G_{g,p}$ 
is a semi-simple algebraic group defined over $\Q$.

\begin{theorem}
\begin{enumerate}
\item If $g\geq 2$, $\rho_p(P\G_g^n)$ and hence $\G_g^n(p)$, is infinite for all odd $p\geq 5$ (\cite{F}).  
\item If $g\geq 2$ and $p\geq 5$ is  either a prime or twice a prime then  $\rho_p$ is irreducible (\cite{Roberts}). Moreover,  $\rho_{p,(i_1,i_2,\ldots,i_m)}$ is irreducible if some $i_s=1$ and $p$ is even (see \cite{KS2}). 
\item $\rho_p$ are asymptotically faithful,  namely $\bigcap_{p_j\to \infty}\ker\rho_{p_j}$ is the centre of the $\G_g^n$ (\cite{A1,FWW}). 
\item The representations $\ro_p\otimes \ro_p^*$  converge in the Fell topology to the action by left composition of $\Gamma_g$ on 
the subspace of regular functions on $\Hom(\pi_1(\Sigma_g), SL(2,\C))$, i.e. the on the vector space generated by the multicurves on $\Sigma_g$ (\cite{CM,MN}). 
\item There exists a free ${\mathcal O}_p$\,-lattice 
$\Lambda_{g,p}$ in $W_g$ and a non-degenerate Hermitian  
${\mathcal O}_p$-valued form on 
$\Lambda_{g,p}$ both invariant under  the action of $\widetilde{\Gamma_g}$ via the representation 
$\widetilde{\rho}_{p,\zeta}$, so that the image of $\ro_{p,\zeta}$ lies in $S\mathbb U_{g,p,\zeta}({\mathcal O}_p)$ when $g\geq 3$  (see \cite{GM}).  
\item $\ro_p(\widetilde{\Gamma_g})$ is a discrete Zariski dense subgroup 
of $\mathbb G_{g,p}(\R)$ whose projections onto the simple factors of 
$\mathbb G_{g,p}(\R)$ are topologically dense, for $g \geq 3$ and  $p\geq 7$ 
prime, $p\equiv 3({\rm mod }\; 4)$ (\cite{LW,F3}). 
\item  If $g\geq 2$ and $p\geq 5$ is prime then  the restriction of $\rho_{p,(2)}$ to $\pi_1(\Sigma_g)\subset \G_g^1$
is infinite (\cite{KS}) and moreover  Zariski dense in $\mathbb G_{g,p}(\R)$ (see \cite{FL}).   
\end{enumerate}
\end{theorem}

The irreducible components of $\rho_p$, for  composite $p$  are not known in full generality, see \cite{AndF,Korinman} for  partial results. 

\begin{remark}\label{mod4}
When $p\equiv 1 ({\rm mod }\; 4)$ the image of the central extension 
of $\Gamma_g$ from \cite{MRo} by $\ro_p$ is contained in 
$\mathbb G_{g,p}(\Z[i])$ and thus it is a discrete 
Zariski dense subgroup of $\mathbb G_{g,p}(\C)$. 
However, if we restrict to the universal central extension  
$\widetilde{\Gamma_g}$  coefficients are reduced from 
$\Z[\zeta_{4p}]$ to $\Z[\zeta_p]$ (see  \cite{GM}, section 13). 
Note that the corresponding invariant form $H_{\zeta_p}$ should be suitably 
rescaled and after rescaling it will be skew-Hermitian when $g$ is odd and Hermitian for even $g$.

\end{remark}

\subsection{Odds and ends}
The AMU Conjecture claims that for a pseudo-Anosov mapping class  $\phi$, i.e. one  which admits a pseudo-Anosov representative homeomorphism,  there exists some $p_0$ such that 
$\rho_p(\phi)$ is of infinite order for every $p\geq p_0$ (see \cite{AMU}). Recall that a homeomorphism of a closed surface is 
pseudo-Anosov if there exist transversely measurable foliations of the surface which are preserved by the homeomorphism while the transverse measures are multiplied by two positive numbers (whose product is 1). Pseudo-Anosov homeomorphisms have precisely two fixed points on the boundary of the Thurston compactification of the Teichm\"uller space.  

Whether the representation of $\G_g^n(p)$ induced by $\rho_p$ is injective seems not known, except for small surfaces (see \cite{FK}). If injectivity holds, then the image of $\rho_p$ cannot be an arithmetic subgroup of $G_{g,p}$ (see \cite{FuPi1}). This would have interesting consequences, as all known infinite images of $\G_g$ are 
arithmetic groups of higher rank, when $g\geq 3$. 
Moreover, it is still unknown whether the representation $\rho_p$ is locally rigid and hence a VHS on a complex projective variety (see \cite{Sim}). We expect the Zariski closure of any MTC representation should be a semisimple Lie group. 
Note that arithmetic lattices of $G_{g,p}$ are cocompact; we don't know whether images of MTC representations can contain unipotent elements. The relation between generalisations of Prym representations (see \cite{Loo}) and quantum representations is still to be uncovered.  

Representations arising for some non semisimple MTC are discussed in \cite{DMJ}, where the authors identified 
them with the homology of local systems on configuration spaces of surfaces with boundary induced from 
Heisenberg group quotients of surface braid groups. This extends the 
construction of Lawrence's representations of braid groups and is a key test for the linearity of mapping class groups.

A version of the volume conjecture (see \cite{CY}), originally due to Kashaev in the case of knots (\cite{Kashaev}),  
states that the logarithmic growth of the Turaev-Viro  invariant at a root of unity $A=\exp(2i\pi/p)$ (for which $H_A$ is not unitary) is the simplicial volume of the closed 3-manifold rescaled by $2\pi$. If true, there would be only finitely many 
hyperbolic 3-manifold with the same skein Turaev-Viro  invariants. Nevertheless there are arbitrarily large finite collections  
of non homeomorphic closed hyperbolic 3-manifolds whose skein Turaev-Viro invariants, volumes, eta and Chern-Simons invariants 
agree (see \cite{KB,Kawauchi,Lickorish2,MNo}). 

\subsection{Outlook}
This survey contains a brief introduction to TQFTs and universal pairings. We discussed manifold detection and positivity of universal pairings, which are bound to work accurately only in low dimensions. We outline the definitions of  fusion, ribbon and  modular tensor categories, which represent the algebraic machinery needed for the construction of TQFTs in dimension 3. We presented some arithmetic properties of modular representations in genus one, in particular the congruence property and discussed partial extensions to 
higher genera, mostly for the skein TQFT. TQFTs in dimension 3 provide a large supply of surface mapping class group representations.  
We expect that  further study of TQFTs can provide new insights about the algebraic properties of mapping class groups.


{
\small      
      
\bibliographystyle{harvard}

\begin{thebibliography}{30}      



\bibitem{A1}
J. E. Andersen, {\em  Asymptotic faithfulness of the 
quantum ${\rm SU}(n)$ representations of the mapping class groups},   
Ann.  Math. (2)  163(2006),  347--368. 

\bibitem{AMU}
J. E. Andersen, G. Masbaum and K. Ueno, {\em 
Topological Quantum Field Theory and the Nielsen–Thurston classification of $M(0,4)$},
Math. Proc. Cambridge Phil. Soc. 141 (2006),  477--488. 

\bibitem{AndF}
J. E. Andersen and J. Fjelstad {\em Reducibility of quantum representations of mapping class groups}, Lett. Math.
Phys. 91 (2010), 215--239, arXiv:0806.2539.



\bibitem{AF}
J. Aramayona and L. Funar 
{\em Quotients of the mapping class group by power subgroups}, 
   Bull. London Math. Soc. 51(2019), 385-398. 

\bibitem{AS}
J. Aramayona and J. Souto, 
{\em Rigidity phenomena in mapping class groups},  Handbook of Teichm\"uller Theory VI 
(A.Papadopoulos, Ed.), IRMA Lect. Math. Theor. Phys. vol. 27,
European Math. Soc., Z\"urich,  2016, 131--165. 




\bibitem{Atiyah}
M. Atiyah, {\em 
Topological quantum field theories}, 
Publ. Math., Inst. Hautes Étud. Sci. 68 (1988), 175--186.


\bibitem{BK}
B.Bakalov and A.Kirillov, Jr.,  
Lectures on tensor categories and modular functors, 
{\em University Lecture Series}, 21. 
American Mathematical Society, Providence, 2001,  x+221 pp. 



\bibitem{B}
P.Bantay, {\em The kernel of the modular representation and the Galois action in RCFT},   
Comm. Math. Phys.  233 (2003),  423–-438.





\bibitem{Barden}
D. Barden, 
{\em Simply connected five-manifolds}, 
Ann. Math. (2) 82 (1965), 365--385. 


\bibitem{BHMV}
C. Blanchet,  N. Habegger, G. Masbaum and P. Vogel, 
{\em  Topological quantum field theories derived from the Kauffman bracket},  
Topology  34 (1995),  883--927.



\bibitem{Br}
M. R. Bridson, {\em Semisimple actions of mapping class groups on ${\rm CAT}(0)$  spaces},  Geometry of Riemann surfaces,  1–14, London Math. Soc. Lecture Note Ser., 368, Cambridge Univ. Press, Cambridge, 2010

\bibitem{CFW}
D. Calegari, M. Freedman and K. Walker, {\em Positivity of the universal pairing in 3 dimensions}, 
J. Amer. Math. Soc. 23 (2009), 109-188. 

\bibitem{CM}
L. Charles and J. March\'e, {\em 
Multicurves and regular functions on the representation variety of a surface in $SU(2)$},  
Commentarii Math. Helv. 87 (2012), 409--431.


\bibitem{CY}
Qingtao Chen, Tian Yang, {\em 
Volume conjectures for the Reshetikhin-Turaev and the Turaev-Viro invariants}, 
Quantum Topology 9 (2018), 419--460.



\bibitem{CG1}
A.Coste and T.Gannon, {\em 
Remarks on Galois symmetry in rational conformal field theories }, Physics Letters B 323 (1994) 316--321.

\bibitem{CG}
A.Coste and T.Gannon, {\em 
Congruence subgroups and rational conformal field theory}, arXiv:math/9909080. 


\bibitem{dBG}
J. De Boer and J. Goeree, {\em Markov traces and II1 factors in conformal field theory}, 
Comm. Math. Phys. 139 (1991),  267--304. 



\bibitem{DM}
P. Deligne and D. Mumford, {\em The irreducibility of the space of curves of given genus}, Publ. Math. I. H. E. S. 36 (1969), 75--109. 

\bibitem{DMo}
P. Deligne and G. D. Mostow, {\em 
Monodromy of hypergeometric functions and non-lattice integral monodromy}, 
Publ. Math. I. H. E. S. 63 (1986), 5--89.


\bibitem{DGGPR}
M. De Renzi, A. Gainutdinov, N. Geer, B. Patureau-Mirand, I. Runkel, {\em 3-Dimensional
TQFTs from Non-Semisimple Modular Categories}, Selecta Math. (N.S.) 28 (2022), no. 2, 42.

\bibitem{DGGPR2}
M. De Renzi, A. Gainutdinov, N. Geer, B. Patureau-Mirand, I. Runkel, {\em Mapping Class Group Representations From Non-Semisimple TQFTs}, Commun. Contemp. Math. 25 (2023), no. 1, 2150091.

\bibitem{DMJ}
M. De Renzi and J. Martel, {\em 
Homological Construction of Quantum Representations of Mapping Class Groups}, arXiv:2212.10940. 



\bibitem{DMa}
B. Deroin and J. March\'e,  
{\em Toledo invariants of Topological Quantum Field Theories}, arXiv:2207.09952








\bibitem{DLN}
Chongying Dong, Xingjun Lin, Xingjun and Siu-Hung Ng, 
{\em Congruence property in conformal field theory}, 
Algebra Number Theory 9 (2015), 2121--2166.

\bibitem{E}
P. Etingof, 
{\em On Vafa’s theorem for tensor categories}, 
Math. Res. Lett. 9 (2002), 651--657.



\bibitem{ENO}
P. Etingof, D. Nikshych and V. Ostrik, 
{\em On fusion categories},
Ann. Math. (2) 162 (2005),  581--642.

\bibitem{ERW}
P. Etingof, E. Rowell and S. Whiterspoon, {\em 
Braid group representations from twisted quantum doubles of finite groups}, 
Pacific J. Math. 234 (2008), 33--41. 


\bibitem{EGNO}
P.Etingof, S. Gelaki, D. Nikshych and V. Ostrik, 
Tensor Categories, MSM 205, AMS 2010. 


\bibitem{EF}
P. Eyssidieux and L. Funar, {\em Orbifold K\"ahler groups related to mapping class groups}, 2020. 


\bibitem{Farb}
B.Farb, Editor, Problems on mapping class groups and related topics. Providence, RI: American Mathematical Society, Proceedings of Symposia in Pure Mathematics 74, 2006.



\bibitem{FaMa}
B. Farb and D. Margalit, {\em A primer of mapping class groups}, Princeton Univ. Press 49, 
2012. 



\bibitem{FF}
J. Fjelstad and J. Fuchs,
{\em Mapping class group representations from Drinfeld doubles of finite groups}, 
J. Knot Theory Ramifications 29, No. 5, Article ID 2050033, 61 p. (2020).






\bibitem{FWW}
M. H. Freedman, K. Walker and Zhenghan Wang, 
{\em Quantum $\rm SU(2)$ faithfully detects mapping class groups modulo center},  Geom. Topol.  6(2002), 523--539. 

\bibitem{FKNSWZ}
M. Freedman, A. Kitaev, C. Nayak, J. Slingerland, K. Walker and Z. Wang, {\em Universal manifold
pairings and positivity}, Geom. Topol. 9 (2005) 2303--2317.

\bibitem{F95}
L. Funar, {\em  2+1-D Topological Quantum Field Theory and 2-D Conformal
Field Theory}, Commun. Math. Phys. 171 (1995), 405--458.

\bibitem{F}
L. Funar, {\em On the TQFT representations of the mapping class groups}, 
Pacific J. Math. 188(1999), 251--274. 


\bibitem{F13}
L. Funar, 
{\em Torus bundles not distinguished by TQFT invariants}, with an Appendix joint with Andrei Rapinchuk,
   Geom. Topol. 17 (2013), 2289--2344.

\bibitem{FK}
L. Funar and T. Kohno, {\em  On Burau representations at roots of unity}, 
   Geom. Dedicata 169(2014), 145--163.

\bibitem{F3}
L. Funar, 
{\em Zariski density and finite quotients of mapping class groups},
   International Mathematics Research Notices 2013, no.9, 2078-2096

\bibitem{FL}
L. Funar and P. Lochak,  {\em Profinite completions of Burnside-type surface groups}, 
   Commun. Math. Phys. 360(2018), 1061--1082. 


\bibitem{FuPi1}
L. Funar and W. Pitsch, {\em 
Images of quantum representations of mapping class groups and Dupont-Guichardet-Wigner quasi-homomorphisms},    J. Inst. Math. Jussieu 17 (2018), 277--304.
   
   
   



\bibitem{FSS}
J. Fuchs, C. Schweigert, C. Stigner, {\em Higher Genus Mapping Class Group Invariants
From Factorizable Hopf Algebras}, Adv. Math. 250 (2014), 285--319. 


\bibitem{Gervais}
 S. Gervais, {\em Presentation and central extensions of mapping class groups}, 
Trans. Amer. Math. Soc. 348 (1996),
3097--3132.










\bibitem{GM}
P. Gilmer and G. Masbaum, {\em Integral lattices in TQFT}, Ann. Sci. E.N.S. 40(2007), 815--844. 














 

\bibitem{Gustafson}
P. Gustafson, 
{\em Finiteness for mapping class group representations from twisted Dijkgraaf–Witten theory}, 
J. Knot Theory Its Ramif. 27 (2018), no. 06, 1850043. 



\bibitem{He}
M. Hennings, {\em Invariants of Links and 3-Manifolds Obtained from Hopf Algebras}, J.
London Math. Soc. (2) 54 (1996),  594--624.


\bibitem{Jones}
V. F. R. Jones, {\em Hecke Algebra Representations of Braid Groups and Link Polynomials}, 
 Ann. Math. 126 (1987), 335--388. 

\bibitem{Juhasz}
A. Juhász, 
{\em Defining and classifying TQFTs via surgery}, 
Quantum Topol. 9 (2018), 229--321.

\bibitem{KB}
J. Kania-Bartoszynska, 
{\em Examples of different 3-manifolds with the same invariants of Witten and Reshetikhin-Turaev}, 
Topology 32 (1993), 47--54.

\bibitem{Kashaev}
R. M. Kashaev, {\em 
The hyperbolic volume of knots from the quantum dilogarithm}, 
Lett. Math. Phys. 39 (1997), 269--275.


\bibitem{KR}
L. H. Kauffman and D. E. Radford, 
{\em Invariants of 3-manifolds derived from finite dimensional Hopf algebras}, 
J. Knot Theory Ramifications 4 (1995), 131--162.

\bibitem{Kawauchi}
A. Kawauchi, {\em 
Topological imitation, mutation and the quantum SU(2) invariants},
J. Knot Theory Ramifications 3 (1994), 25--39.


\bibitem{KL}
T. Kerler, V. Lyubashenko, Non-Semisimple Topological Quantum Field Theories
for 3-Manifolds with Corners, Lecture Notes in Mathematics 1765. Springer-Verlag,
Berlin, 2001.



\bibitem{KS}
T. Koberda and R. Santharoubane,
{\em Quotients of surface groups and homology of finite covers via quantum representations}, 
Invent. Math. 206 (2016), 262--292. 


\bibitem{KS2}
T. Koberda and R. Santharoubane,
{\em
Irreducibility of quantum representations of mapping class groups with boundary}, 
Quantum Topol. 9 (2018), 633--641. 


\bibitem{Korinman}
J. Korinman, {\em 
Decomposition of some Reshetikhin-Turaev representations into irreducible factors}, 
SIGMA 15 (2019), 011, 25 pages.

\bibitem{KT}
M. Kreck and P. Teichner, {\em Positivity of topological field theories in dimension at least 5}, J. Topology 1 (2008), 663--670. 



\bibitem{Kup}
G. Kuperberg,
{\em  Non-involutory Hopf algebras and 3-manifold invariants}, Duke Math. J. 84 (1996), 83--129. 



\bibitem{LW}
M. Larsen and  Zhenghan Wang, {\em  
Density of the SO(3) TQFT representation of mapping class groups},  
Comm. Math. Phys.  260(2005),  641--658. 





\bibitem{Lickorish}
W. B. R.Lickorish, {\em 
Invariants for 3-manifolds from the combinatorics of the Jones polynomial}, 
Pac. J. Math. 149 (1991),  337--347.

\bibitem{Lickorish4}
W. B. R.Lickorish, {\em 
Skeins and handlebodies}, 
Pac. J. Math. 159 (1993), 337--349.

\bibitem{Lickorish2}
W. B. R.Lickorish, {\em 
Distinct 3-manifolds with all $SU(2)_q$ invariants the same}, 
Proc. Am. Math. Soc. 117 (1993), 285--292.


\bibitem{Lickorish3}
W. B. R.Lickorish, {\em 
What is … a skein module?}, 
Notices Am. Math. Soc. 56 (2009), 240--242.


\bibitem{Liu}
Yi Liu, {\em 
Finite-volume hyperbolic 3-manifolds are almost determined by their finite quotient groups}, Invent. Math. to appear, 
arXiv:2011.09412. 






\bibitem{Loo}
E. Looijenga,  {\em Prym representations of mapping class groups}, Geom. Dedicata 64 (1997), 69--83. 




\bibitem{Ly}
V. Lyubashenko, {\em Invariants of 3-Manifolds and Projective Representations of Mapping
Class Groups via Quantum Groups at Roots of Unity}, Comm. Math. Phys. 172 (1995), 467--516. 

\bibitem{Ly2}
V. Lyubashenko,
{\em Ribbon abelian categories as modular categories}, 
J. Knot Theory Ramifications 5 (1996),  311--403.


\bibitem{MN}
J. March\'e and  M. Narimannejad, 
{\em Some asymptotics of topological quantum field theory via skein theory}, 
Duke Math. J. 141 (2008),  573--587. 


\bibitem{Marche}
J. March\'e, Introduction to quantum representations of mapping class groups, notes 
de cours, Bordeaux, 2018. 

\bibitem{MRo}
G. Masbaum and J. D. Roberts, {\em 
On central extensions of mapping class groups},  
Math. Ann.  302(1995), 131--150.

\bibitem{MNo}
S. V. Matveev and T. Nowik, {\em 
On 3-manifolds having the same Turaev-Viro invariants}
Russ. J. Math. Phys. 2 (1994), 317--324.


\bibitem{Mcmullen}
C. T. McMullen, {\em Braid groups and Hodge theory}, Math. Annalen 355 (2013), 893--946. 



\bibitem{MS}
G. Moore and N. Seiberg, 
{\em  Classical and quantum conformal field theory}, Commun.
Math. Phys. 123 (1989), 177--254.






\bibitem {NS}
Siu-Hung Ng and P.Schauenburg, {\em Congruence subgroups and generalized 
Frobenius-Schur indicators}, Comm. Math. Phys. 
300 (2010), 1--46.








\bibitem{RT}
N. Reshetikhin and V. G. Turaev, {\em Invariants of 3-manifolds via link polynomials and quantum groups}, Invent. Math. 103 (1991),  547--597.


\bibitem{Roberts1}
J. Roberts, {\em 
Skeins and mapping class groups}, 
Math. Proc. Camb. Philos. Soc. 115 (1994), 53--77.


\bibitem{Roberts}
J. Roberts, {\em Irreducibility of some quantum representations of mapping class groups}, 
Knots in Hellas '98, Vol. 3 (Delphi). 
J. Knot Theory Ramifications 10 (2001), no. 5, 763--767. 






























\bibitem{Segal}
G. B. Segal, {\em  The definition of conformal field theory}, Differential geometrical methods in theoretical physics, Proc. 16th Int. Conf., NATO Adv. Res. Workshop, Como/Italy 1987, NATO ASI Ser., Ser. C 250, 165--171, 1988.


\bibitem{Segal2}
G. B. Segal, 
{\em The definition of conformal field theory}, 
Tillmann, Ulrike (ed.), Topology, geometry and quantum field theory. Proceedings of the 2002 Oxford symposium in honour of the 60th birthday of Graeme Segal, Oxford, UK, June 24–29, 2002. Cambridge: Cambridge University Press (ISBN 0-521-54049-6/pbk). London Mathematical Society Lecture Note Series 308, 421--577 (2004).

\bibitem{Sim}
C. T. Simpson, {\em Higgs bundles and local systems}, Publ. Math. I. H. E. S. 
 75 (1992), 5--95. 





\bibitem{Turaev}
V. Turaev, Quantum invariants of knots and 3-manifolds, De Gruyter Studies in Mathematics 18, 
2016. 

 




\bibitem{Witten}
E. Witten, {\em Quantum field theory and the Jones polynomial}, Comm. Math. Physics 121 (1989),  351--399.





\bibitem{Xu}
Feng Xu, {\em Some computations in the cyclic permutations of completely rational nets},   Comm. Math. Phys.  267  (2006),   757--782.




\end{thebibliography}

}

\end{document}